\documentclass[12pt]{article}

\usepackage{myarticle}
\usepackage[colorinlistoftodos,prependcaption,textsize=scriptsize]{todonotes}

\usepackage{mycommands}
\usepackage{bm}

\colorlet{changecolor}{orange}

\newcommand{\opt}[1]{#1^*}

\newcommand{\ak}[1]{\bar{#1}}
\newcommand{\akk}[1]{\tilde{#1}}
\newcommand{\akp}[1]{\hat{#1}}

\newcommand{\piter}[1]{_{#1}}

\newcommand{\bmat}[1]{{\bm #1}}
\newcommand{\spdmat}{\mathbb S_{++}}
\newcommand{\spsdmat}{\mathbb S_{+}}

\newcommand{\Tau}{T}
\newcommand{\lbstep}{a}

\newcommand{\nbeta}{R}

\DeclareMathOperator*{\dist}{dist}

\newcommand{\matid}{{\bm\opid}}

\newcommand{\matM}{{\bm M}}
\newcommand{\matQ}{{\bm Q}}

\newcommand{\DM}{\bmat D_{\bmat M}}

\renewcommand{\rto}[1]{\overset{#1}{\to}}

\newcommand{\lin}[4]{\ell_{#1}^{#2}(#3,#4)}
\newcommand{\linfun}[2]{\ell_{#1}^{#2}}    

\newcommand{\KL}{Kurdyka--{\L}ojasiewicz\xspace}

\newif\ifShowFigures
\ShowFigurestrue

\graphicspath{{./figures/}}

\begin{document}

\thispagestyle{empty}
\begin{center}
\vspace*{0.03\paperheight}
  {\Large\bf
    Adaptive FISTA for Non-convex Optimization
  }\\ \bigskip \bigskip
  {\large 
    Peter Ochs${}^*$ and Thomas Pock${}^\dagger$
  } \\ \medskip
  {\small
    $^*$~Saarland University, Saarbr\"{u}cken, Germany \\
    $^\dagger$~TU Graz, Graz, Austria
  }
\end{center}
\bigskip

\begin{abstract}
  In this paper we propose an adaptively extrapolated proximal gradient method, which is based on the accelerated proximal gradient method (also known as FISTA), however we locally optimize the extrapolation parameter by carrying out an exact (or inexact) line search. It turns out that in some situations, the proposed algorithm is equivalent to a class of SR1 (identity minus rank 1) proximal quasi-Newton methods. Convergence is proved in a general non-convex setting, and hence, as a byproduct, we also obtain new convergence guarantees for proximal quasi-Newton methods. The efficiency of the new method is shown in numerical experiments on a sparsity regularized non-linear inverse problem.
\end{abstract}
\bigskip


\section{Introduction}

The introduction of accelerated gradient methods by Nesterov in~\cite{Nest83} has arguably revolutionized the world of large-scale convex and non-smooth optimization. Computationally, they are as simple as plain gradient descent, but come along with a much faster rate of convergence. Therefore, they found numerous applications in modern signal/image processing, and machine learning. It is well-known that these accelerated gradient methods are \enquote{optimal} for the class of smooth and convex (not necessarily strongly-convex) problems with Lipschitz continuous gradient, in the sense that their worst case complexity is proportional to the theoretical lower complexity bound of first-order methods for this class of problems \cite{NY83,Nest04}. The mechanism of accelerated gradient methods can be interpreted and explained in many different ways and hence, the magic of acceleration is still the subject of intensive research efforts. We refer to \cite{Goh17}, for a detailed exposition of several aspect of acceleration.

Accelerated gradient methods can be seen as a modification of the heavy-ball method of Polyak~\cite{Polyak64}. While the heavy-ball method already achieves an optimal convergence rate on smooth and strongly convex functions with Lipschitz continuous gradient, accelerated gradient methods are optimal, among first-order methods, for smooth convex problems with Lipschitz continuous gradient~\cite{Nest04}.

In the original work~\cite{Nest04}, the acceleration is explained by the concept of so called estimation sequences, which generate a sequence of simple convex functions approximating the original function. Accelerated gradient methods also show strong connections to the St\"ormer--Verlet method~\cite{Hairer03} for discretizing second order ordinary differential equations (ODEs). Using this connection, both the heavy-ball method~\cite{Polyak64} and Nesterov's original method~\cite{Nest83} can be seen as particular types of discrete-time approximations to the \enquote{heavy-ball with friction} dynamical system. The main principle here is to attach a \enquote{mass} to the sequence of points generated which accelerates when moving down the landscape of the objective function. More recent works use this relation to investigate the properties of accelerated gradient methods in the framework of ODEs~\cite{Su14,Attouch16}. Finally, accelerated first-order methods can also be explained by accelerated dual schemes~\cite{Nest07} or by accelerated primal-dual algorithms~\cite{CP11}.

For quadratic functions, both the heavy-ball method and accelerated
gradient methods show striking similarities to the conjugate gradient
method~\cite{Stiefel1952}. However, while in the conjugate gradient
method, the step size and extrapolation parameters are chosen locally
using an exact line (or plane) search, accelerated gradient methods
derive their parameters from global properties of the objective such
as the smallest and largest eigenvalues of the corresponding system
matrix. On the one hand, the conjugate gradient method is optimal, among first-order methods, and
comes along with a finite convergence property for quadratic functions.  On the other hand,
its generalization to non-quadratic functions is much harder. At
this point, accelerated gradient methods are advantageous. They
can be applied not only to quadratic problems but to the whole class
of differentiable functions with Lipschitz continuous gradient.

Accelerated gradient methods have been generalized to be applicable to
a class of convex optimization problems that can be written as the sum
of a differentiable function with Lipschitz continuous gradient and a
non-smooth function with easy to compute proximal map. The most
popular instance is the Fast Iterative Shrinkage Thresholding
Algorithm (FISTA)~\cite{BT09b} but also other, more general schemes
have been proposed. See the work of Tseng~\cite{Tseng08} for an
excellent presentation and unification of a whole family of
accelerated proximal gradient methods. 

Finally, accelerated gradient methods have also turned out to perform
very well on non-convex optimization problems. However, while their
empirical performance is often comparable to the convex case, their
convergence properties are still hardly understood from a theoretical
point of view~\cite{GL15,LL15,DK16,OCBP14}.

In this paper, we study a modification of FISTA for
minimizing an objective given by the sum of a smooth function $f$ with
Lipschitz continuous gradient and a non-smooth function $g$ with
simple proximal mapping. We still rely on a step size parameter
related to the global Lipschitz parameter $L$, however similar to the
conjugate gradient method, we propose to locally adapt the
extrapolation parameter by means of an exact line search. A simplified
version of the algorithm's update step for minimizing a smooth function $f$ is the following:
\[
  \begin{split}
    y\idx{\beta}\iter\k =&\  x\iter\k + \beta (x\iter\k - x\iter\km) \\
    x\iter\kp\idx{\beta} \in&\ \argmin_{x\in \R^N}\, g(x) 
    + \scal{\nabla f(y\idx{\beta}\iter\k)}{x-y\idx{\beta}\iter\k} 
    + \frac L2 \vnorm{x - y\idx{\beta}\iter\k}^2 \,,
  \end{split}
\]
where $x_k$ is the main iteration, $y_k^{(\beta)}$ is the extrapolated
point (depending on the choice of $\beta)$ and $\beta \in \R$ is a
free parameter which is optimized in each step in order to provide the
locally largest decrease of the proximal subproblem. It turns
out that if the function $f$ is quadratic, the optimal choice of
$\beta$ makes the scheme equivalent to a proximal SR1 (identity minus
rank one) quasi-Newton method. This equivalence yields a connection
between accelerated schemes and quasi-Newton methods, which holds even
when $g$ is a non-convex function. The convergence that we establish
for our method in a general non-smooth non-convex setting translates
directly to a class of proximal quasi-Newton methods. Unlike the
general class of proximal quasi-Newton methods which usually suffer
from the problem that simple proximal mappings become hard to solve,
efficient solutions to the \enquote{rank-1 proximal mappings} are
known \cite{BF12,KV17,BFO18}.\\

The remainder of the paper is organized as follows: In Section~\ref{sec:rel-work} we provide a review of state-of-the-art and put the contributions of this paper into related work. In Section~\ref{sec:opt-problem}, we present the considered optimization problem, and in Section~\ref{sec:alg}, we detail the proposed algorithm and give different convergence guarantees for both convex and non-convex problems. Numerical results of the proposed algorithm are provided in Section~\ref{sec:numerics}. In the last section we give some conclusions and discuss open problems for future research.

\section{Related work} \label{sec:rel-work}

\paragraph{Proximal Gradient Method.} The basic update step of our method is a so-called proximal gradient step (also known as forward-backward splitting) \cite{LM79,CW05,CP11a,Nest04}. Although, using this basic step yields an efficient algorithm in many cases, it has been observed that the worst case complexity is not optimal \cite{Nest04} for certain classes of convex optimization problems. This observation led to the exploration of so-called accelerated schemes. The basic proximal gradient step can be accelerated by an additional (computationally cheap) extrapolation step \cite{BT09b,Nest83,Tseng08} (see also \cite{Tseng10}), where the iteration-dependent extrapolation parameter obeys a certain rule derived from global properties of the objective function. In contrast to the prescribed rules, in this paper, we explicitly optimize the extrapolation parameter. Surprisingly, for a certain class of problems, the optimized scheme has a closed form expression as a proximal quasi-Newton method, which yields new convergence results.

\paragraph{Classical quasi-Newton Methods.} Quasi-Newton methods are intensively studied in the classic context of non-linear optimization problems. We refer to \cite{NW06,DS96b,Broyden67} for an overview and references. The basic idea of quasi-Newton methods is to successively improve a quadratic approximation to the objective function, i.e., the goal is to approximate second order information (Hessian) using combinations of first order information. Maybe the most widely known and used quasi-Newton method is BFGS \cite{Fletcher87} or its limited memory variant L-BFGS \cite{LN89}, and its extension to bound constraints L-BFGS-B \cite{BLNZ94}. More results that motivate the acceleration governed by such a variable metric approach are \cite{BGLS95,BQ99,PLS08}.

\paragraph{Quasi-Newton Methods for non-smooth problems.} Though, originally designed for smooth optimization problems, the BFGS method shows good performance on non-smooth problems as well. In \cite{YVGN10}, the method was reinterpreted and analysed in the non-smooth convex setting. Theoretical guarantees of the original BFGS method for one or two dimensional non-smooth problems are established in \cite{LO13,LZ15}. Motivated by the original global convergence proof of Powell \cite{Powell76}, Guo and Lewis \cite{GL17} substantially extended the previous theoretical guarantees to a class of non-smooth convex problems. In general, convergence cannot be expected.

\paragraph{Proximal quasi-Newton Methods.} Driven by the success of proximal splitting methods, the concept of quasi-Newton methods was also applied to optimization problems with more structure than just smoothness, similar to the setting of Proximal Gradient Descent. The crucial aspect for the efficiency of such a variable metric Proximal Gradient Method is the evaluation of the proximal mapping, which is often expensive when coordinates are not separated in the objective.  There are specific choices of the metric that allow the proximal mapping to be solved efficiently, e.g. when the metric is a diagonal matrix, which preserves the separability. Becker and Fadili \cite{BF12} derived efficient solutions when the metric is of type \enquote{identity plus rank 1}.  The proximal step was embedded into a zero memory variant of the SR1 quasi-Newton method \cite{NW06}.  Unlike the BFGS-method, which updates the Hessian approximation with a matrix of rank 2, the SR1 method updates the approximation with a rank 1 matrix.  Nevertheless, Nocedal and Wright \cite{NW06} state the observation that the \enquote{[...] SR1 method appears to be competitive with the BFGS method}. A special case of an \enquote{identity minus rank 1} metric was studied in \cite{KV17}, which also leads to an efficiently solvable proximal mapping. Recently, proximal mappings with respect to metrics of type ``identity plus/minus rank $\nbeta$'' have been studied in~\cite{BFO18}. In \cite{FG16}, interior point methods are used to solve the proximal mapping efficiently for so-called quadratic support functions.\\

Variable metric versions of the Proximal Gradient Method have been studied without paying special attention to efficiently solving the proximal mapping. The earliest reference is \cite[Section 5]{CR97}. In the broader context of monotone inclusion problems, the general framework for analyzing the convergence of variable metric methods is that of quasi-Fej\'{e}r sequences \cite{CV13}, which was used in \cite{CV14} to find a relation to primal--dual schemes. For a variable metric algorithm with mild differentiability assumption (without the usual gradient Lipschitz continuity), we refer to \cite{Salzo16}.  The convergence of forward--backward splitting with iteration dependent Bregman distances in an infinite dimensional setting, which contains the variable metric proximal step as a special case, was proved in \cite{Nguyen17}.

\paragraph{Non-convex setting.} The variable metric Proximal Gradient Method in \cite{CPR13} follows the majorization--minimization principle \cite{HL04}. The metric is constructed to induce a quadratic majorizer in each iteration. However, the arising proximal mappings are computationally expensive. Extensions to a block-coordinate descent version are presented in \cite{CPR16} and a combination with an inertial method in \cite{Ochs16}. In \cite{BLPP16,BLPPR16}, the choice of the metric enjoys a great flexibility at the cost of an additional line search step in the algorithm. An extension of their framework to a general non-smooth first-order oracle and a flexible choice of Bregman distances was proposed in \cite{OFB17}. A different way to incorporate a variable metric into forward--backward splitting (FBS) was proposed in \cite{STP17}. They reinterpret FBS as a gradient descent method with a variable metric, which allows them to use some of the machinery from smooth optimization. Besides accelerating optimization methods by a variable metric, accelerated (optimal) gradient methods can be used in the non-convex setting \cite{GL15,LL15,DK16} as well, where the goal is to guarantee the accelerated rate of convergence when the same method is applied to convex problems.

\paragraph{Subspace approach.} The proposed method shows similarities to sequential subspace optimization (SESOP) \cite{NZ05,EMZ07} (see also the review article \cite{ZE10}) and the Majorize-Minimize (MM) subspace algorithm \cite{CIM11,CJPT13} (when $g\equiv0$), which have been introduced as effective acceleration strategies. They generate the next iterate by minimizing the objective (if computationally affordable) or a constructed surrogate function over a subspace. Many variants to define these subspace have been considered. In contrast, we optimize for the extrapolation point on a subspace from where a gradient descent step yields the best decrease of a simple quadratic surrogate function. In general, both approaches differ, however majorizers constructed in subspace approaches could be incorporated also into our scheme.

\subsection{Contribution}

\paragraph{Convergence of an optimally extrapolated Proximal Gradient Method (PGM).} We study a novel, optimally extrapolated PGM for non-convex optimization problems. If the function is composed of a continuously differentiable function with Lipschitz continuous gradient and a non-smooth function, we prove that a subsequence converges to a stationary point and that objective values converge under reasonable conditions. 

\paragraph{Equivalence to proximal quasi-Newton methods.} In a broad special non-convex setting, we prove our method to be equivalent to a proximal variant of the SR1 quasi-Newton method \cite{NW06}. To be more precise, we prove the equivalence between an optimally rank-$\nbeta$ extrapolated Proximal Gradient Method and an identity minus rank-$\nbeta$ proximal quasi-Newton method. Our convergence results can be applied directly to these methods, leading to the first convergence result of this SR1 proximal quasi-Newton for non-convex optimization problems.

\paragraph{Relation to \cite{BF12} and \cite{KV17}.} A related proximal quasi-Newton method was considered by Becker and Fadili in \cite{BF12}. However, they adapt a Barzilai--Borwein step length, whereas we relate the step length to the Lipschitz constant. Moreover, the proposed quasi-Newton metric in the evaluation of the proximal mapping is of type \enquote{identity plus a rank 1 matrix}, whereas ours is of type \enquote{identity minus rank 1}, and we allow for non-convex functions in the proximal mapping. Another closely related approach is that of Karimi and Vavasis \cite{KV17}, which applies to problems that are the sum of a convex quadratic and a convex non-smooth function. The considered proximal mapping is also of type \enquote{identity minus rank 1}, however the relation to the SR1 metric \cite[Section 8.2]{NW06} is unclear and the analysis is completely different to ours. The case of an \enquote{identity plus or minus a rank-$\nbeta$ matrix}, which generalizes \cite{BF12,KV17}, was recently studied in \cite{BFO18}.


\section{The Optimization Problem} \label{sec:opt-problem}

The setting of this paper is that of a Euclidean vector space $\R^N$ of dimension $N$ equipped with the standard inner product $\scal\cdot\cdot$ and induced norm $\vnorm{x} = \sqrt{\scal xx}$ for $x\in \R^N$. Moreover, we use the notation $\vnorm[\bmat V]{x}^2:= \scal{x}{x}_{\bmat V} := \scal{x}{\bmat V x}$ for a matrix $\bmat V\in \R^{N\times N}$. We denote by $\spsdmat(N)$ the set of symmetric positive semi-definite matrices in $\R^{N\times N}$ and by $\spdmat(N)\subset\spsdmat(N)$ the subset of positive definite matrices.\\

We consider the optimization problem
\begin{equation} \label{eq:gen-opt}  \tag{P}
  \min_{x\in \R^N}\, f^g(x)\,,\quad f^g(x):= g(x) + f(x) 
\end{equation}
with the following assumptions.
\begin{ASS} \label{ass:problem}
Let $\bmat L\in\spdmat(N)$. 
\begin{enumerate}
  \item\label{ass:problem:i} $\map{g}{\R^N}{\eR}$ is proper lower semi-continuous (lsc).
  \item\label{ass:problem:ii} $\map{f}{\R^N}{\R}$ is continuously differentiable and $\nabla (f\circ \bmat L^{-1/2})$ is $1$-Lipschitz continuous.
  \item\label{ass:problem:iii} $f^g$ is bounded from below.
\end{enumerate}
\end{ASS}
\begin{REM}
  Assumption~\ref{ass:problem}(\ref{ass:problem:ii}) is equivalent to $1$-Lipschitz continuity of the gradient of $f$ in a metric induced by $\bmat L$ (see Lemma~\ref{lem:Lipschitz-grad-equiv}), and implies a Generalized Descent Lemma (see Lemma~\ref{lem:gen-Descent-Lemma}).  
\end{REM}

\section{The Adaptive Extrapolation Algorithm aFISTA}

Before we introduce the full algorithm in Section~\ref{sec:alg}, we present the main update step in Section~\ref{sec:update-step}. In a special setting, the update step is shown to be equivalent to a certain proximal quasi-Newton method without extrapolation in Section~\ref{sec:rel-to-quasi-Newton}. This relation reveals many interesting application, for which the main update step of our algorithm can be computed efficiently. Section~\ref{sec:mm} describes an approach that allows more general problems to be tackled with the same strategy. Additional strategies for the efficient computation of the update are discussed in Section~\ref{sec:efficient-update}. Convergence of the main algorithm is analyzed in Section~\ref{sec:convergence-analysis}. All proofs are given in the appendix.

\subsection{The Extrapolated Proximal Gradient Update Step} \label{sec:update-step}

Let $\ak{x}\in\R^N$ be the current point of the iterative scheme and $\set{d_1,\ldots,d_\nbeta}$ be vectors in $\R^N$ and
\begin{equation} \label{eq:extrapolate-y}
  y\idx{\bvec\beta}:= \ak{x} + \sum_{i=1}^\nbeta \beta_i d_i
\end{equation}
for $\bvec{\beta}:=(\beta_1,\ldots,\beta_\nbeta) \in \R^\nbeta$. In matrix--vector notation, we have $y\idx{\bvec\beta}= \ak{x} + \bmat D\bvec \beta$ where $\bmat D=(d_1,\ldots,d_\nbeta)\in \R^{N\times \nbeta}$ contains the vectors $d_1,\ldots,d_\nbeta$ as columns. We denote the by
\begin{equation} \label{eq:def-g-linearization}
    \ell_f^g(x; y\idx{\bvec\beta}) := g(x) + f(y\idx{\bvec\beta}) + \scal{\nabla f(y\idx{\bvec\beta})}{x-y\idx{\bvec\beta}},
\end{equation}
the linearized smooth function $f$ plus the non-smooth function
$g$. The extrapolated proximal gradient (EPG) step is given by the
following combined $x$ and $\bvec\beta$ proximal problem:
\begin{equation} \label{eq:aFista-step-gen-ho}   \tag{EPG}
    (\akp{x},\akp{\bvec\beta})\in \argmin_{x\in \R^N} \min_{\bvec\beta\in \R^\nbeta}\, \ell_f^g(x; y\idx{\bvec\beta}) + \frac 1{2} \vnorm[\bmat\Tau]{x-y\idx{\bvec\beta}}^2 \,,
\end{equation}
where $\bmat\Tau\in\spdmat(N)$. We
also allow for inexact minimizers $(\akk{x}, \akk{\bvec\beta})$ of
\eqref{eq:aFista-step-gen-ho} as long as the following
condition is satisfied:
\begin{equation} \label{eq:aFista-inexact-step-gen-ho}  \tag{iEPG}
    \ell_f^g(\akk{x}; y\idx{\akk{\bvec\beta}}) + \frac 1{2} \vnorm[\bmat\Tau]{\akk{x}-y\idx{\akk{\bvec\beta}}}^2 \leq f^g(\ak{x}).
\end{equation}
Observe that the right hand side of \eqref{eq:aFista-inexact-step-gen-ho} is the value of the objective in \eqref{eq:aFista-step-gen-ho} at $\ak{x}$ with $\beta=0$. Obviously, \eqref{eq:aFista-step-gen-ho} implies \eqref{eq:aFista-inexact-step-gen-ho}, i.e., $\akk{x}=\akp{x}$ and $\akk{\bvec\beta}=\akp{\bvec\beta}$ with $(\akp{x},\akp{\bvec\beta})$ from \eqref{eq:aFista-step-gen-ho} asserts the inexact condition in \eqref{eq:aFista-inexact-step-gen-ho}.

\subsection{The Algorithm} \label{sec:alg}

We embed the update step from Section~\ref{sec:update-step} into our main algorithm. While we state Algorithm~\ref{alg:afista-nc} with the exact update step \eqref{eq:aFista-step-gen-ho}, the convergence analysis in Section~\ref{sec:convergence-analysis} also accounts for the inexact version \eqref{eq:aFista-inexact-step-gen-ho}. 
\mybox{
  \begin{ALG}[Adaptive Fista for Non-convex Problems with Backtracking] \ \label{alg:afista-nc}
    \vspace*{-3ex}
    \begin{itemize}
      \item \emph{Optimization Problem}: \eqref{eq:gen-opt} with Assumption~\ref{ass:problem}.
      \item \emph{Initialization}: $x\iter0=\bar x$ for some $\bar x\in \R^N$. Set $\lbstep>0$.
      \item \emph{Selection Rule}: Select extrapolation directions $\seq[\k\in\N]{\bmat D\iter\k}$ with $\bmat D\iter\k\in \R^{N\times \nbeta}$. 
      \item \emph{Update for $\k\geq 0$}: Find $x\iter\kp$ and $\bvec\beta\iter\k$ such that the following holds:
      \begin{equation} \label{eq:update:alg:afista-nc} \tag{EPG}
            (x\iter\kp,\bvec\beta\iter\k) \in \argmin_{x\in\R^N}\, \min_{\bvec\beta\in \R^\nbeta}\, \ell_f^g(x; y\idx{\bvec\beta}) + \frac 1{2} \vnorm[\bmat\Tau\piter\k]{x-y\idx{\bvec\beta}}^2 
      \end{equation}
      with $y\idx{\bvec\beta} = x\iter\k + \bmat D\iter\k \bvec\beta$, where $\bmat\Tau\piter\k - \bmat L\piter\k - \lbstep\matid \in\spsdmat(N)$ and $\bmat L\piter\k$ satisfies
      \begin{equation} \label{eq:update:alg:afista-nc-L-cond}
        f(x\iter\kp) \leq f(y\idx{\bvec\beta\iter\k}) + \scal{\nabla f(y\idx{\bvec\beta\iter\k})}{x\iter\kp - y\idx{\bvec\beta\iter\k}} + \frac 12 \vnorm[\bmat L\piter\k]{x\iter\kp - y\idx{\bvec\beta\iter\k}}^2 \,.
      \end{equation}
    \end{itemize}
  \end{ALG}
}
\begin{REM}
  The extrapolation directions $\bmat D\iter\k$ may, for example, depend on $\seq[j\leq\k]{x\iter{j}}$ or $\k$.
\end{REM}
The underlying scheme of the algorithm is the extrapolated proximal gradient algorithm. It is complemented with a locally varying Lipschitz upper bound that resembles the Descent Lemma. However, we also allow for an anisotropic quadratic upper bound, represented by the metric $\bmat\Tau\piter\k$. This may yield more flexibility and tighter bounds.
\begin{EX}
  Let $\bmat\Tau\piter\k=\bmat L\piter\k = L\cdot\matid$ with identity matrix $\matid$ where $L>0$ is the Lipschitz constant of $\nabla f$ and $\bmat D\iter\k=x\iter\k-x\iter\km$. Then, \eqref{eq:update:alg:afista-nc-L-cond} is satisfied and, for $\bvec\beta\iter\k=\beta\piter\k\in\R$ satisfying\footnote{Since \eqref{eq:update:alg:afista-nc} minimizes over $x$ and $\bvec\beta$, the optimal $\bvec\beta$ may actually depend on $x$.} \eqref{eq:update:alg:afista-nc}, the update step \eqref{eq:update:alg:afista-nc} can be written as follows:
  \[
    \begin{split}
    y\idx{\beta\piter\k} =&\ x\iter\k + \beta\piter\k(x\iter\k - x\iter\km) \\ 
    x\iter\kp \in &\  \argmin_{x\in \R^N}\, g(x) + \frac L{2} 
                 \vnorm{x-\big(y\idx{\beta\piter\k} - \frac 1L \nabla f(y\idx{\beta\piter\k})\big)}^2 \,,
    \end{split}
  \]
  which is an extrapolated proximal gradient step. In the case $f$ and $g$ are convex functions, and $\beta\piter\k$ is suboptimal, depending on $k$  like $\seq[\k\in\N]{1-3/\k}$, formally, we obtain the accelerated proximal gradient method, which is known as FISTA \cite{BT09b}. It achieves a convergence rate of $O(1/\k^2)$ in term of objective value residual.

As the proposed update step in \eqref{eq:update:alg:afista-nc} takes the same form, but optimizes the update with respect to $\bvec\beta$, each step guarantees a better objective value than an accelerated proximal gradient step. However, this is not enough to guarantee the same rate of convergence, which requires a global picture of the objective.
\end{EX}

\subsection{Relation to Proximal Quasi-Newton Methods} \label{sec:rel-to-quasi-Newton}

For a quadratic function $f$, the optimal parameter $\bvec \beta$ can be computed analytically. In Section~\ref{sec:mm}, we describe a strategy that applies this result in a more general setting. Using the optimal $\bvec\beta$ in the update step yields an interpretation of the update step \eqref{eq:aFista-step-gen-ho} as a proximal gradient step in a modified metric without extrapolation, also known as proximal quasi-Newton method. The modified metric is of type ``identity minus rank $\nbeta$'', where $\nbeta$ is the column rank of $\bmat D$. Proximal quasi-Newton methods have recently drawn attention, as for $\nbeta=1$, the associated proximal mapping can be evaluated efficiently \cite{BF12,KV17}. See~\cite{BFO18}, for the rank-$\nbeta$ case.
\begin{THM} \label{thm:aFista-equiv-ho}
  Let Assumption~\ref{ass:problem} hold and consider the problem in \eqref{eq:aFista-step-gen-ho}. Suppose that $f(x) = \frac 12 \scal{x}{\bmat H x}+\scal{b}{x} + c$ is quadratic with Hessian $\bmat H\in \spdmat(N)$, $b\in \R^N$, $c\in \R$, the matrix $\bmat\Tau$ is chosen such that $\bmat M:=\bmat\Tau -\bmat H \in \spdmat(N)$, and the columns of $\bmat D\in\R^{N\times \nbeta}$ are linearly independent.  Then, the inner optimization problem w.r.t. $\bvec \beta$ is solved by
\[
  \bvec\beta^* = (\bmat D^\top \bmat M \bmat D)^{-1} \bmat D^\top \bmat M (x-\ak{x})\,.
\]
The optimization problem in \eqref{eq:aFista-step-gen-ho} is equivalent to the following:
\begin{equation} \label{eq:thm:aFista-equiv-ho}
  \akp{x} \in\argmin_{x\in \R^N}\, g(x) + \frac 12 \vnorm[\bmat Q]{x- \ak{x} + \bmat Q^{-1}\nabla f(\ak{x})}^2 \,,
\end{equation}
where 
\begin{equation} \label{eq:thm:aFista-equiv-ho-Qmatrix}
  \matQ:= \bmat\Tau  - \bmat U^\top\bmat U \in \spdmat(N) \qquad \text{with}\qquad \bmat U := (\bmat D^\top \bmat M \bmat D)^{-\frac 12} \bmat D^\top \bmat M
\end{equation}
and $\bmat U^\top \bmat U$ is of rank $\nbeta$. The inverse metric is
\begin{equation} \label{eq:thm:aFista-equiv-ho-Qmatrix-inv}
  \bmat Q^{-1} = \bmat\Tau ^{-1} + \bmat\Tau ^{-1} \bmat U^\top (\matid - \bmat U \bmat\Tau ^{-1} \bmat U^\top)^{-1} \bmat U \bmat\Tau ^{-1} \,.
\end{equation}
\end{THM}
The proof is in Section~\ref{proof:thm:aFista-equiv-ho}.
\begin{REM}
  The update step \eqref{eq:thm:aFista-equiv-ho} can be written equivalently as
  \begin{equation} \label{eq2:thm:aFista-equiv-ho}
    \akp{x} \in\argmin_{x\in\R^N}\, g(x) + f(\ak{x}) + \scal{\nabla f(\ak{x})}{x- \ak{x}} +  \frac 12 \vnorm[\bmat Q]{x- \ak{x}}^2 \,.
  \end{equation}
\end{REM}
\begin{REM} \label{rem:interpret-Q}
  The proof of Theorem~\ref{thm:aFista-equiv-ho} shows that, on $\im(\bmat D)$, which denotes the image of $\bmat D$, the matrix $\bmat Q$ coincides exactly with the Hessian matrix, i.e., $\bmat Q x = \bmat H x$ for $x\in \im(\bmat D)$, or equivalently, $\bmat Q\bmat D \bvec\beta=\bmat H\bmat D\bvec\beta$ for $\bvec\beta\in\R^\nbeta$. On the orthogonal complement $\im(\bmat D)^\bot=\ker(\bmat D^\top)$, the matrix $\bmat Q$ coincides with $\bmat T$, i.e. $\bmat Q x = \bmat T x$ for $x$ with $\bmat D^\top x=0$. 
\end{REM}
\begin{REM} \label{rem:know-H-on-imQ}
  The formulas for the matrix $\bmat Q$ in \eqref{eq:thm:aFista-equiv-ho-Qmatrix} and $\bmat Q^{-1}$ in \eqref{eq:thm:aFista-equiv-ho-Qmatrix-inv} do not require explicit inversion of $\bmat H$. In fact, $\bmat H$ must be known on $\im(\bmat D)$ only. Using $\bmat Y:=\bmat H \bmat D$, the formula in \eqref{eq:thm:aFista-equiv-ho-Qmatrix} can be written as 
  \[
    \bmat Q = \bmat T - (\bmat T \bmat D - \bmat Y) \big[\bmat D^\top (\bmat T\bmat D - \bmat Y) \big]^{-1} (\bmat T\bmat D - \bmat Y)^\top  
  \]
  and \eqref{eq:thm:aFista-equiv-ho-Qmatrix-inv} as 
  \[
    \bmat Q^{-1} = \bmat T^{-1} + (\bmat D - \bmat T^{-1} \bmat Y)\big[(\bmat D - \bmat T^{-1} \bmat Y)^\top \bmat Y\big]^{-1}(\bmat D - \bmat T^{-1} \bmat Y)^\top \,.
  \]
  The details of the computation are given in Section~\ref{proof:rem:know-H-on-imQ}. These formulas provide an elegant and efficient way to compute that application of $\bmat Q$ and its inverse.
\end{REM}

\paragraph{The rank 1 case.} In practical applications, the case where $\bvec\beta=\beta$ is 1-dimensional is most interesting, since the resulting quasi-Newton method is of type ``identity minus rank 1'',\linebreak for which the update steps can be evaluated efficiently. In fact, it turns out that in this case \eqref{eq:aFista-step-gen-ho} is equivalent to a generalization of the SR1 quasi Newton method with a non-smooth term, i.e., a proximal SR1 quasi Newton method. 
\begin{COR} \label{cor:aFista-equiv}
Consider the situation of Theorem~\ref{thm:aFista-equiv-ho} with $\nbeta=1$, i.e. $\bmat D = d\in \R^N$. Then, the inner optimization problem w.r.t. $\bvec\beta=\beta\in\R$ is solved by
\[
  \opt\beta = \frac{\scal{d}{x-\ak{x}}_\matM}{\scal{d}{d}_\matM}\,.
\]
The optimization problem in \eqref{eq:aFista-step-gen-ho} is equivalent to the following:
\begin{equation} \label{eq:id-minus-rank1-fbs-step}
  \akp{x} \in\argmin_{x\in \R^N}\, g(x) + \frac 12 \vnorm[\bmat Q]{x- \ak{x} + \matQ^{-1}\nabla f(\ak{x})}^2 \,,
\end{equation}
where 
\begin{equation} \label{eq:id-minus-rank1-metric}
  \bmat Q:= \bmat\Tau  - u u^\top \quad \text{with}\quad u := \frac{\matM d}{\vnorm[\matM]{d}} 
\quad\text{and}\quad 
  \bmat Q^{-1} = \bmat\Tau ^{-1} + \frac{\bmat\Tau ^{-1} u u^\top \bmat\Tau ^{-1}}{1-u^\top \bmat\Tau ^{-1} u} \,. 
\end{equation}
\end{COR}
\begin{proof}
  The statement is an obvious consequence of Theorem~\ref{thm:aFista-equiv-ho}.
\end{proof}
\begin{REM} \label{rem:prox-SR1-metric}
  We explain the naming \enquote{proximal SR1 quasi-Newton method}. The standard notation of quasi-Newton methods uses $\overline{\bmat B}\piter\k$ for the current approximation of the Hessian matrix and $\overline{\bmat H}\piter\k$ for the approximation of the inverse Hessian matrix. We use the ``overline'' to distinguish $\overline{\bmat H}\piter\k$ from the second derivative $\bmat H$ of $f$ in Theorem~\ref{thm:aFista-equiv-ho} and Corollary~\ref{cor:aFista-equiv}. The SR1 update formula \cite[Section 8.2]{NW06} is given by
\[
  \begin{split}
  \overline{\bmat B}\piter\kp 
  =&\  \overline{\bmat B}\piter\k + \frac{(y\iter\k-\overline{\bmat B}\piter\k d\iter\k)(y\iter\k-\overline{\bmat B}\piter\k d\iter\k)^\top}{\scal{d\iter\k}{y\iter\k-\overline{\bmat B}\piter\k d\iter\k}}\,, \\
  \overline{\bmat H}\piter\kp 
  =&\  \overline{\bmat H}\piter\k + \frac{(d\iter\k-\overline{\bmat H}\piter\k y\iter\k)(d\iter\k - \overline{\bmat H}\piter\k y\iter\k)^\top}{\scal{d\iter\k-\overline{\bmat H}\piter\k y\iter\k}{y\iter\k}} \,. 
  \end{split}
\]
  Setting $\nbeta=1$, $d=x\iter\k-x\iter\km$, $y=\bmat H d=\nabla f(x\iter\k) - \nabla f(x\iter\km)$ in \eqref{eq:id-minus-rank1-metric} yields:
\[
  \bmat Q = \bmat\Tau  + \frac{(y-\bmat\Tau  d)(y-\bmat\Tau d)^\top}{\scal{d}{y-\bmat\Tau  d}}\quad\text{and}\quad
  \bmat Q^{-1} = \bmat\Tau ^{-1} + \frac{(d-\bmat\Tau ^{-1} y)(d - \bmat\Tau ^{-1}y)^\top}{\scal{d-\bmat\Tau ^{-1} y}{y}} \,,
\]
which is that same as Remark~\ref{rem:know-H-on-imQ} specified for $\nbeta =1$. The identification $\overline{\bmat B}\piter\k = \bmat T$ and $\overline{\bmat B}\piter\kp=\bmat Q$ shows that aFISTA in \eqref{eq:id-minus-rank1-fbs-step} uses a $0$-memory SR-1 metric, i.e., we re-initialize $\overline{\bmat B}\piter\k$ in each iteration with $\bmat T$.

  Unfortunately, the memory based SR1 update formula does not preserve the positive definiteness, which makes the convergence analysis for line search based algorithms difficult. For more details on this discussion, we refer to \cite[Section 8.2]{NW06}. Our $0$-memory variant guarantees positive definiteness of $\bmat Q$.

  In case we restart the approximation of the Hessian matrix in each iteration with a multiple of the identity matrix (with a Barzilai--Borwein rule), we obtain the recently proposed zero memory version of the SR1 quasi-Newton method \cite{BF12}. 
\end{REM}

\subsection{A Majorization--Minimization Strategy} \label{sec:mm}

Using the majorization--minimization strategy \cite{HL04}, we can apply the results from the previous section to more general problems. This idea is also pursued in \cite{CIM11,CJPT13} for smooth (possibly non-convex) functions (i.e. $g\equiv 0$). However, we consider the full non-convex setting (see Assumption~\ref{ass:problem}). At the current iterate $\ak{x}$, we define a quadratic tangent majorizer\footnote{A tangent majorizer $Q_f(x,\ak{x})$ coincides with $f$ at $x=\ak{x}$, has the same slope as $f$ at $\ak{x}$ and $Q_f(x,\ak{x})\geq f(x)$ holds for all $x$.} of $f$, for example, 
\begin{equation} \label{eq:mm:our-quad-tangent-majorizer}
  Q^{\bmat L}_f (x, \ak{x}) := f(\ak{x}) + \scal{\nabla f(\ak{x})}{x-\ak{x}} + \frac 12 \vnorm[\bmat L]{x-\ak{x}}^2\,,
\end{equation}
which is chosen according to the quadratic upper bound from the Generalized Descent Lemma (Lemma~\ref{lem:gen-Descent-Lemma}). Let $\akp{x}$ be the solution of one iteration of \eqref{eq:aFista-step-gen-ho} applied to $g(x) + Q^{\bmat L}_f (x, \ak{x})$ with $\bmat T$ such that $\bmat T-\bmat L \in \spdmat(N)$. This setting allows us to apply Theorem~\ref{thm:aFista-equiv-ho} with $\bmat Q:= \bmat T - \bmat U^\top \bmat U$ and $\bmat H=\bmat L$ to compute $\akp{x}$. Then, the following relation holds:
\[
  \begin{split}
  f^g(\akp{x}) 
  \overset{\ii1}{\leq}&\  g(\akp{x}) + Q^{\bmat L}_f(\akp{x}, \ak{x}) \\
  \overset{\ii2}{\leq}&\  g(\akp{x}) + Q^{\bmat Q}_f(\akp{x},\ak{x}) - \frac{1}{2}\vnorm[\bmat Q-\bmat L]{\akp{x}-\ak{x}}^2 \\
  \overset{\ii3}{\leq}&\  g(\ak{x}) + Q^{\bmat Q}_f(\ak{x},\ak{x}) - \frac{1}{2}\vnorm[\bmat Q-\bmat L]{\akp{x}-\ak{x}}^2 \\
  \leq&\  f^g(\ak{x}) - \frac{1}{2}\vnorm[\bmat Q-\bmat L]{\akp{x}-\ak{x}}^2 \\
  \end{split}
\]
where $\ii1$ uses the properties of a tangent majorizer, $\ii2$ uses Theorem~\ref{thm:aFista-equiv-ho}, and $\ii3$ uses \eqref{eq2:thm:aFista-equiv-ho}. Remark~\ref{rem:interpret-Q} shows that 
\[
    \vnorm[\bmat Q-\bmat L]{\akp{x}-\ak{x}}^2 = \vnorm[\bmat T - \bmat L]{\akp{x}^0-\ak{x}^0}^2 \geq c \vnorm{\akp{x}^0-\ak{x}^0}^2
\]
holds for some $c>0$, where $\akp{x}^0-\ak{x}^0$ is the projection of $\akp{x}-\ak{x}$ onto $\ker(\bmat D^\top)$, which is a sufficient decrease condition. If a sufficient decrease condition with respect to all coordinates of $\akp{x}-\ak{x}$ is desired, we can perform the update step \eqref{eq:thm:aFista-equiv-ho} with $\bmat Q = \bmat T - \rho \bmat U^\top \bmat U$ for some $\rho\in [0,1)$.

The algorithm is summarized in Algorithm~\ref{alg:mm-afista-nc}.
\mybox{
  \begin{ALG}[MM based Adaptive Fista for Non-convex Problems] \ \label{alg:mm-afista-nc}
    \begin{itemize}
      \item \emph{Optimization Problem}: \eqref{eq:gen-opt} with Assumption~\ref{ass:problem}.
      \item \emph{Initialization}: $x\iter0=\bar x$ for some $\bar x\in \R^N$. Set $\rho\in[0,1)$.
      \item \emph{Selection Rule}: Select extrapolation directions $\seq[\k\in\N]{\bmat D\iter\k}$ with $\bmat D\iter\k\in \R^{N\times \nbeta}$. 
      \item \emph{Update for $\k\geq 0$}: Define a tangent majorizer $Q^{\bmat L\iter\k}_f(x,x\iter\k)$ around $x\iter\k$, \linebreak set $\bmat Q\iter\k=\bmat T\iter\k - \rho \bmat U\iter\k^\top \bmat U\iter\k$ with $\bmat U\iter\k$ as in Theorem~\ref{thm:aFista-equiv-ho}, and compute $x\iter\kp$ as follows:
      \begin{equation} \label{eq:update:alg:afista-nc-MM}
        x\iter\kp \in  \argmin_{x\in \R^N}\, g(x) + \frac 12 \vnorm[\bmat Q\iter\k]{x- x\iter\k + \bmat Q\iter\k^{-1}\nabla f(x\iter\k)}^2
      \end{equation}
      where $\bmat T\iter\k - \bmat L\iter\k \in \spdmat(N)$.
    \end{itemize}
  \end{ALG}
}
The convergence of this algorithm was analyzed in the context of variable metric forward--backward splitting in \cite{CPR13} when $g$ is a convex function. Assuming the \KL inequality holds, which is basically satisfied by any function appearing in practical applications, they show convergence of the full sequence generated by the algorithm to a stationary point. \cite{CPR16} covers also the case of a non-convex function $g$. 
\begin{REM}
  In practice, there are several ways to construct tangent majorizers. This issue has been considered, for example, in \cite{CIM11,CJPT13} in the context of subspace optimization algorithm. Another possibility is given by approximating the Hessian matrix at $x\iter\k$ by a symmetric positive definite matrix, which can be achieved, for example, by adding a multiple of the identity matrix or by projection onto $\eps\cdot\matid+\spdmat(N)$ for some $\eps>0$, if the Hessian is not positive definite \cite{NW06}.
\end{REM}

\subsection{Efficient Computation of the Update Step} \label{sec:efficient-update}

Here, we discuss the efficient computation of the update step used in Algorithms~\ref{alg:afista-nc} and~\ref{alg:mm-afista-nc}. 
\begin{itemize}
  \item If $f$ is a convex quadratic function or Algorithm~\ref{alg:mm-afista-nc} is used, and $\bmat T\iter\k$ is a diagonal matrix, the update step is a proximal step of type ``identity minus rank $\nbeta$''. According to \cite{BFO18}, such proximal mappings can be computed efficiently. The computation can be reduced to proximal mappings with respect to the diagonal part of the metric and $\nbeta$-dimensional root finding problems, which can be solved efficiently thanks to their small dimension.
  \item In the same setting as above, but with $\nbeta=1$, \cite{BF12} provides a list of efficiently computable rank-1 proximal mapping, which includes $g$ being the $\ell_1$-norm, group $\ell_{1,2}$-norms, $\ell_\infty$-norm, simplex, box, or non-negativity constraints, etc. However, note that we also allow for non-convex functions $g$. 
  \item Unless we can find analytic formulas in special cases, the update step in \eqref{eq:update:alg:afista-nc} needs to be computed numerically. For example, this can be achieved by alternating minimization. The minimization with respect to $\bvec\beta$ is a $\nbeta$-dimensional smooth minimization problem, and the minimization with respect to $x$ is a proximal gradient step around the extrapolated point $y\idx{\bvec\beta}$. The subproblem w.r.t. $\bvec \beta$ is convex (quadratic), if $f$ in \eqref{eq:update:alg:afista-nc} is convex (quadratic). In particular, this is true when the majorization--minimization strategy (Section~\ref{sec:mm}) is used.
  \item Finally, we can also allow to inexactly computed update steps. The main requirement is a sufficient improvement on the objective function in \eqref{eq:update:alg:afista-nc} from $x\iter\k$ to $x\iter\kp$. This allows for an efficient backtracking strategy for $\bvec\beta$, which always succeeds with $\bvec\beta=0$. Therefore, we can try a few samples of $\bvec\beta$ before falling back to $\bvec\beta=0$, e.g. for $\nbeta=1$, we can use Nesterov's sequence $\beta\piter\k\to1$ for $\k\to\infty$.
\end{itemize}

\subsection{Convergence Analysis} \label{sec:convergence-analysis}

Let us analyze Algorithms~\ref{alg:afista-nc} and~\ref{alg:mm-afista-nc}. Preliminaries from non-smooth analysis are provided in Section~\ref{appdx:prelim}.

\begin{PROP}[Convergence of Adaptive Fista] \label{prop:conv-nc}
  Let $\seq[\k\in\N]{x\iter\k}$, $\seq[\k\in\N]{\bvec\beta\iter\k}$, $\seq[\k\in\N]{\bmat D\iter\k}$, $\seq[\k\in\N]{\bmat L\iter\k}$ and $\seq[\k\in\N]{\bmat T\iter\k}$ be as in Algorithm~\ref{alg:afista-nc} for solving \eqref{eq:gen-opt} with Assumption~\ref{ass:problem}.  Suppose $\seq[\k\in\N]{\bmat T\piter\k}$ is bounded. Then, the sequence of objective values $\seq[\k\in\N]{f^g(x\iter\k)}$ is non-increasing and every limit point of $\seq[\k\in\N]{x\iter\kp}$ is a stationary point of $f^g$.
\end{PROP}
The proof is in Section~\ref{proof:prop:conv-nc}.
\begin{REM}
  The boundedness assumption for $\seq[\k\in\N]{\bmat T\piter\k}$ can always be satisfied thanks to Assumption~\ref{ass:problem}(\ref{ass:problem:ii}).
\end{REM}
The existence of a convergent subsequence converging to a stationary point can be derived analogously to Proposition~\ref{prop:conv-nc}. As we have mentioned above, the convergence of such variable metric proximal gradient algorithm is known \cite{CPR13}.

In the following statement, we prove convergence for an inexact version of Algorithm~\ref{alg:afista-nc}. Of course, this more general setting requires a few more assumptions.
\begin{PROP}[Convergence of Inexact Adaptive Fista] \label{prop:conv-inc}
  Let $\seq[\k\in\N]{x\iter\k}$, $\seq[\k\in\N]{\bvec\beta\iter\k}$, $\seq[\k\in\N]{\bmat D\iter\k}$, $\seq[\k\in\N]{\bmat L\iter\k}$ and $\seq[\k\in\N]{\bmat T\iter\k}$ be as in Algorithm~\ref{alg:afista-nc} with \eqref{eq:update:alg:afista-nc} replaced by \eqref{eq:aFista-inexact-step-gen-ho}, i.e., find $x\iter\kp$ and $\bvec\beta\piter\k$ that satisfy
  \begin{equation} \label{prop:conv-inc:inexact-update}
    \ell_f^g(x\iter\kp; y\idx{\bvec\beta\piter\k}) + \frac 1{2} \vnorm[\bmat\Tau]{x\iter\kp-y\idx{\bvec\beta\piter\k}}^2 \leq f^g(x\iter\k)\,,
  \end{equation}
  for solving \eqref{eq:gen-opt} with Assumption~\ref{ass:problem}.  Suppose $\seq[\k\in\N]{\bmat T\piter\k}$ is bounded. Then, the sequence of objective values $\seq[\k\in\N]{f^g(x\iter\k)}$ is non-increasing. If, additionally, $g$ is continuous on its domain and $\norm[-]{\partial f^g(x\iter\kp)}\to 0$ for $\k\to\infty$, then every limit point of $\seq[\k\in\N]{x\iter\kp}$ is a stationary point of $f^g$.
\end{PROP}
The proof is in Section~\ref{proof:prop:conv-inc}.

\subsubsection{Discussion of the Convergence Rate}
  In the convex setting, we can easily derive a convergence rate of $O(1/\k)$ w.r.t. the function gaps for Algorithms~\ref{alg:afista-nc} and~\ref{alg:mm-afista-nc}. However, we cannot guarantee an optimal convergence rate of $O(1/\k^2)$ as for FISTA. Unlike FISTA, our algorithm performs a local optimization, and hence does not exploit the global properties of convex functions. Interestingly, the conjugate gradient method also performs a local optimization (in both the step size and the extrapolation parameter) and still guarantees an optimal $O(1/\k^2)$ convergence rate. Hence it will be interesting in future work to study an extension of the proposed algorithm which also performs a locally optimal choice of the step size parameter.

  However, we can derive hybrid accelerated algorithms that make use of the proposed update step \eqref{eq:aFista-step-gen-ho}. Since, the main focus of this paper is the equivalence between the adaptive extrapolation method and quasi-Newton methods, we do not provide further details here. For the sake of completeness, we list the mentioned variants in the appendix in Section~\ref{sec:acc-aFista} and prove their optimal rate of convergence.

\section{Numerical Experiments} \label{sec:numerics}

For optimization problems \eqref{eq:gen-opt} where $f$ is a quadratic function and $g$ is convex, adaptive FISTA (aFISTA) is equivalent to the zero memory SR1 quasi-Newton method with step size $\bmat\Tau\piter\k = \tau\piter\k^{-1} \matid$ and $\tau\piter\k < 1/L$, instead of Barzilai--Borwein step sizes. Therefore, its efficiency has already been demonstrated in \cite{BF12}.  In Section~\ref{sec:exp-nn}, we consider a non-convex and non-smooth setting. 

\subsection{Sparsity regularized non-linear inverse problem}  \label{sec:exp-nn}

In this experiment, we consider a neural network formulation of a one
dimensional regression problem. We are given $N=80$ noisy samples $(X,\tilde Y)\in 
(\R^{1\times N})^2$ of the function $\map{F}{[-3,3]}{\R}$, $x\mapsto x^3
+ \cos(5x)$, arranged as corresponding columns of the matrices $X$ and $\tilde
Y$, i.e., 
\[
  \tilde Y_{1,i} := F(X_{1,i}) + E_{1,i} \,,
\]
where $E\in \R^{1\times N}$ is an additive noise matrix that models Gaussian
noise with standard deviation $\frac 32$ and 20 randomly scaled outliers. 
The neural network optimization problem (non-linear inverse problem)
that we consider is formulated as follows:
\begin{gather} \label{eq:min-prob-SimpleNet}
  \min_{\substack{W_0, W_1, W2\\ b_0, b_1, b_2}}\, \sum_{i=1}^N \Big(\vnorm{(W_2\sigma_2(W_1 \sigma_1(W_0 X+B_0 )+B_1)+B_2 - \tilde Y)_{1,i} }^2 + \eps^2 \Big)^{1/2} + \lambda \sum_{j=0}^2 \norm[1]{W_j} \,, \\
  B_j := b_j \bvec 1^\top\,,\quad \bvec 1^\top := (1,\ldots, 1)\,, \quad j=1,2,3\,,
\end{gather}
where $D_0=D_3=1$, $D_1=D_2=10$, $W_j\in \R^{D_j\times D_{j-1}}$, $B_{j+1}\in
\R^{D_j\times N}$, for $j=1,2,3$, and $\map{\sigma_j}{\R^{D_j\times
N}}{\R^{D_j\times N}}$, $A\mapsto (\max(0,(A_{i,l}^2+\eps^2)^{1/2}))_{i,l}$,
for $j=1,2$, with $\eps = 0.1$. We used $\lambda=1$, which led to a sparsity
level of about $87\%$ of the coordinates.

\begin{figure}[t]
  \begin{center}
\ifShowFigures
    \begin{tikzpicture}
      \begin{axis}[%
        width=0.48\linewidth,%
        height=7.5cm,%
        xmode=log, ymode=log,
        xlabel={iteration $\k$},%
        ylabel=$h(x\iter\k) / h(x\iter0)$,%
        legend columns=5,
        legend cell align=left,
        legend style={at={(1.0,1.02)},anchor=south,font=\footnotesize,column sep=10pt}
        ]
      \addplot[very thick,blue!80!white,solid,mark=triangle,mark repeat={50},mark size={2},mark options={solid}] %
          table {figures/SimpleSparseNet_conv_FBS_iter.dat};
          \addlegendentry{FBS};
      \addplot[very thick,red!80!black,solid,mark=star,mark repeat={50},mark size={2},mark options={solid}] 
          table {figures/SimpleSparseNet_conv_iPiano_iter.dat};
          \addlegendentry{iPiano};
      \addplot[very thick,yellow!80!green,solid,mark=square,mark repeat={50},mark size={2},mark options={solid}] %
          table {figures/SimpleSparseNet_conv_mFISTA_iter.dat};
          \addlegendentry{MFISTA};
      \addplot[very thick,green!80!black,solid,mark=*,mark repeat={50},mark size={2},mark options={solid}] 
          table {figures/SimpleSparseNet_conv_aFISTA_iter.dat};
          \addlegendentry{aFISTA};
      \addplot[very thick,cyan!80!black,solid,mark=o,mark repeat={50},mark size={2},mark options={solid}] 
          table {figures/SimpleSparseNet_conv_zeroSR1_PG_iter.dat};
          \addlegendentry{ZeroSR1};
      \end{axis}
      \begin{scope}[xshift=0.47\linewidth]
      \begin{axis}[%
        width=0.48\linewidth,%
        height=7.5cm,%
        xmode=log, 
        xlabel={time [sec]}%
        ]
      \addplot[very thick,blue!80!white,solid,mark=triangle,mark repeat={50},mark size={2},mark options={solid}] %
          table {figures/SimpleSparseNet_conv_FBS_time.dat};
      \addplot[very thick,red!80!black,solid,mark=star,mark repeat={50},mark size={2},mark options={solid}] 
          table {figures/SimpleSparseNet_conv_iPiano_time.dat};
      \addplot[very thick,yellow!80!green,solid,mark=square,mark repeat={50},mark size={2},mark options={solid}] %
          table {figures/SimpleSparseNet_conv_mFISTA_time.dat};
      \addplot[very thick,green!80!black,solid,mark=*,mark repeat={50},mark size={2},mark options={solid}] 
          table {figures/SimpleSparseNet_conv_aFISTA_time.dat};
      \addplot[very thick,cyan!80!black,solid,mark=o,mark repeat={50},mark size={2},mark options={solid}] 
          table {figures/SimpleSparseNet_conv_zeroSR1_PG_time.dat};
      \end{axis}
            \end{scope}
    \end{tikzpicture} 
\else
  \texttt{Activate Figure Later: Comparison Simple Net}
\fi
  \end{center}
  \caption{\label{fig:Comparison-SimpleNet}Convergence plots for the methods
  described in Section~\ref{sec:exp-nn} for solving
  \eqref{eq:min-prob-SimpleNet}. The vertical axis is the same for both plots.
  In this experiment, our method aFISTA outperforms the other methods with
  respect to the actual computation time and the final objective value.}
\end{figure}
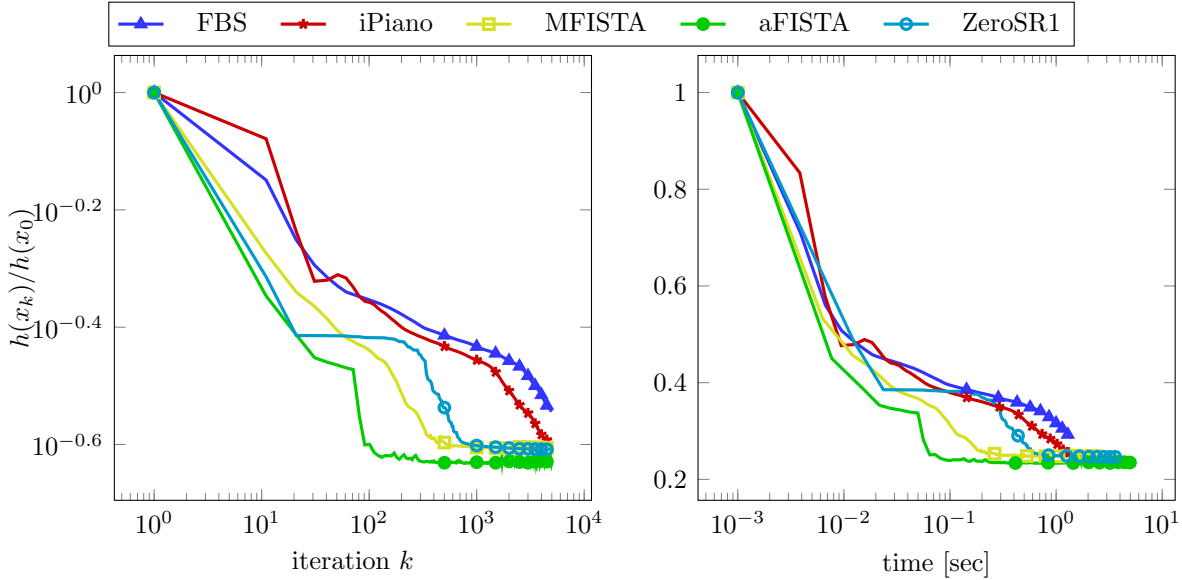

We run Algorithm~\ref{alg:afista-nc}, denoted aFISTA, with inexact update \eqref{eq:aFista-inexact-step-gen-ho} and backtracking strategy for the extrapolation ($\beta\piter\k\in\set{2,1,0}$ and $\bmat D\iter\k=x\iter\k-x\iter\km$), and Algorithm~\ref{alg:mm-afista-nc}, denoted ZeroSR1, with majorizer computed using the (heuristic) SR1 metric (cf. Remark~\ref{rem:prox-SR1-metric}), which leads to the algorithm in \cite{BF12} (with line-search). Besides our methods, we evaluate Forward--Backward Splitting (FBS), iPiano with $\beta=0.95$ \cite{OCBP14}, and monotone FISTA (MFISTA) \cite{LL15}, a variant of FISTA with convergence guarantees in the non-convex setting. We used the same heuristic step size $\alpha=5\cdot 10^{-5}$ for all methods.  Throughout several experiments, the performance of aFISTA was on a par with MFISTA. The convergence for one problem instance is shown in Figure~\ref{fig:Comparison-SimpleNet}. In this experiment, aFISTA is slightly better and finds a lower objective value. Since aFISTA used a small number of backtracking iterations, the performance with respect to time is also competitive with MFISTA. The ZeroSR1 method is competitive in the beginning. 

\section{Conclusion}

In this paper, we analyzed a non-convex variant of the well-known FISTA, where
the extrapolation parameter is adaptively optimized in each iteration, which we
call aFISTA. In the special case where the smooth part of the split objective
function is quadratic, aFISTA is equivalent to a certain proximal
quasi-Newton method, which unlike the general class of proximal quasi-Newton
methods, allows for efficient solutions of the proximal mapping. This
equivalence relies on a reformulation of the quadratic function and is not
influenced by the non-smooth part of the objective, which may also be
non-convex. It provides a different view on quasi-Newton methods, which allows
for accelerated variants. 

The general convergence of aFISTA is studied for non-convex objective functions that are the sum of a continuously differentiable function with Lipschitz continuous gradient and a proper lower semi-continuous non-smooth function with simple proximal mapping. We prove that a subsequence converges to a stationary point in terms of the limiting subdifferential. In numerical experiments, aFISTA and its variants have been shown to be competitive with the state-of-the-art.

In future work, we will explore the relationship between the adaptive 
extrapolation strategy (aFISTA) with $\nbeta$ linearly independent directions and 
other well known quasi-Newton methods such as BFGS. Of course, this 
requires efficient solution algorithms for rank-$\nbeta$ proximal mappings. Moreover,
our perspective of the SR1 quasi-Newton method might also lead to new convergence
results when the construction of the metric uses memory of previous Hessian 
approximations. As a third direction of future research, the usage of Bregman 
distances needs to be investigated. For strongly convex Bregman functions, 
convergence should not be difficult to proof, however, it is not clear which
\enquote{rank-$1$} Bregman proximal mappings can be solved efficiently.

\appendix


\section{Preliminaries}

We present here a version of the Descent Lemma that uses Lipschitz continuity in a specific metric. For equivalent characterizations of the Lipschitz condition, see Lemma~\ref{lem:Lipschitz-grad-equiv}.
\begin{LEM}[Generalized Descent Lemma] \label{lem:gen-Descent-Lemma}
  Let $\bmat V\in \spdmat(N)$ and $\map{f}{\R^N}{\R}$ be a function. Then $\nabla_{\bmat V} f:=\bmat V^{-1}\nabla f$ is $L$-Lipschitz continuous in the metric induced by $\bmat V$, i.e. 
  \begin{equation} \label{eq:Lipschitz-cond-M}
      \norm[\bmat V]{\nabla_{\bmat V}f( x) - \nabla_{\bmat V} f( \bar x)}  \leq  L \norm[\bmat V]{ x -  \bar x}  \quad \forall  x, \bar x\in \R^N 
  \end{equation}
  holds, then the following inequality holds:
  \begin{equation} \label{eq:lem:gen-Descent-Lemma}
      \abs{f(x) - f(\bar x) - \scal{\nabla f(\bar x)}{x-\bar x}} \leq \frac L2\norm[\bmat V]{x-\bar x}^2 \quad \forall  x, \bar x\in \R^N\,.
  \end{equation}
\end{LEM}
\begin{proof}
  Using the Fundamental Theorem of Calculus, Cauchy--Schwarz inequality in the metric $\bmat V$, and the Lipschitz continuity of $\nabla_{\bmat V}f$, we obtain the result:
  \[
    \begin{split}
      \abs{f(x) - f(\bar x) - \scal{\nabla f(\bar x)}{x-\bar x}}
      =&\ \abs{\int_0^1 {\scal{\bmat V^{-1}(\nabla f(\bar x + t(x-\bar x)) - \nabla f(x))}{x-\bar x}}_{\bmat V}\, \mathit{d}t \,} \\
      \leq&\ \int_0^1 {\norm[\bmat V]{\nabla_{\bmat V} f(\bar x + t(x-\bar x)) - \nabla_{\bmat V} f(x)}\norm[\bmat V]{x-\bar x}}\, \mathit{d}t \\
      \leq&\ \int_0^1 {tL\norm[\bmat V]{x-\bar x}^2}\, \mathit{d}t = \frac L2\norm[\bmat V]{x-\bar x}^2 \,.
      \qedhere
    \end{split}
  \]
\end{proof}
\begin{LEM} \label{lem:Lipschitz-grad-equiv}
  Let $\bmat V\in \spdmat(N)$ and $\map{f}{\R^N}{\R}$ and set $g:=f\circ \bmat V^{-1/2}$. The Lipschitz continuity condition in \eqref{eq:Lipschitz-cond-M} is equivalent to each of the following statements:
  \begin{itemize}
    \item $L$-Lipschitz continuity of $\nabla g$
    \item $\nabla f$ is $L$-Lipschitz w.r.t. $(\R^N,\vnorm[\bmat V]{\cdot})$ and (the dual space\footnote{The dual norm of $\norm[\bmat V]{\cdot}$ is $\norm[\bmat V]{\cdot}^*=\vnorm[\bmat V^{-1}]{\cdot}$.}) $((\R^N)^*, \vnorm[\bmat V^{-1}]{\cdot})$, i.e.
    \[
        \norm[\bmat V^{-1}]{\nabla f(x)-\nabla f(\bar x)} \leq L\norm[\bmat V]{x-\bar x}\quad \forall  x, \bar x\in \R^N\,.
    \]
  \end{itemize}
\end{LEM}
\begin{proof}
  The result follows from the following equivalences:
  \[
    \begin{array}{rrcll}
      & \norm{\nabla g (x) - \nabla g(y)} & \leq & L \norm{x - y} &\quad \forall x,y\in \R^N  \\
      \Leftrightarrow 
      & \norm{\bmat V^{-1/2} (\nabla f(\bmat V^{-1/2}x) - \nabla f( \bmat V^{-1/2}y))}&  \leq&  L \norm{x - y} & \quad \forall x,y\in \R^N  \\
      \Leftrightarrow 
      & \norm{\bmat V^{-1/2} (\nabla f(\tilde x) - \nabla f(\tilde y))}&  \leq&  L \norm{\bmat V^{1/2}(\tilde x - \tilde y)} & \quad \forall \tilde x,\tilde y\in \R^N  \\
      \Leftrightarrow 
      & \norm{\bmat V^{1/2}(\bmat V^{-1}\nabla f(\tilde x) - \bmat V^{-1}\nabla f(\tilde y))}&  \leq&  L \norm{\bmat V^{1/2}(\tilde x - \tilde y)} & \quad \forall \tilde x,\tilde y\in \R^N  \\
      \Leftrightarrow 
      & \norm[\bmat V]{\bmat V^{-1}\nabla f(\tilde x) - \bmat V^{-1}\nabla f(\tilde y)}&  \leq&  L \norm[\bmat V]{\tilde x - \tilde y} & \quad \forall \tilde x,\tilde y\in \R^N \\
      \Leftrightarrow 
      & \norm[\bmat V]{\nabla_{\bmat V}f(\tilde x) - \nabla_{\bmat V} f(\tilde y)}&  \leq&  L \norm[\bmat V]{\tilde x - \tilde y} & \quad \forall \tilde x,\tilde y\in \R^N 
    \end{array}
  \]
\end{proof}

\subsection{Concepts from non-smooth analysis} \label{appdx:prelim}

We summarize preliminaries from \cite{Rock98}. The \emph{Fr\'echet subdifferential} of $\map{f}{\R^N}{\eR}$ at $\bar x \in\dom f:=\set{x\in\R^N\,\vert\, f(x)<+\infty}$ is the set $\rpartial f(\bar x)$ of those elements $v \in \R^N$ such that
\[
  \lim\inf_{\substack{x\to \bar x\\ x\neq \bar x}} \frac{f(x) - f(\bar x) - \scal{v}{x-\bar x}}{\vnorm{x-\bar x}} \geq 0  \,.
\]
For $\bar x\not\in \dom f$, we set $\rpartial f(\bar x) = \emptyset$. The \emph{(limiting) subdifferential} of $f$ at $\bar x\in\dom f$ is defined by
\[
  \partial f(\bar x) := \set{v\in \R^N\vert\, \exists\, (x\iter\k,f(x\iter\k)) \to (\bar x, f(\bar x)),\;v\iter\k\in \rpartial f(x\iter\k),\;v\iter\k \to v} \,,
\]
and $\partial f(\bar x) = \emptyset$ for $\bar x \not\in \dom f$. A point $\bar x\in \dom f$ for which $0\in \partial f(\bar x)$ is a called a \emph{critical points}. As a direct consequence, we have the following closedness property: 
\[
  (x\iter\k,f(x\iter\k)) \to (\bar x, f(\bar x)),\ v\iter\k\to \bar v,\ \text{and for all } \k\in\N\colon v\iter\k \in \partial f(x\iter\k)\quad \Longrightarrow\quad  \bar v\in \partial f(\bar x) \,.
\]
The \emph{distance} of $\bar x\in\R^N$ to a set $\omega\subset \R^N$ as is given by $\dist(\bar x,\omega) := \inf_{x\in\omega}\, \vnorm{\bar x - x}$. We introduce the \emph{non-smooth slope} $\vnorm[-]{\partial f(\bar x)}:= \inf_{v\in\partial f(\bar x)} \vnorm{v} =\dist(0, \partial f(\bar x))$ at $\bar x$. Note that $\inf \emptyset := +\infty$ by definition. Furthermore, we have (see \cite{FGP14}):
\begin{LEM} \label{lem:lazy-slope-liminf}
If $(x\iter\k,f(x\iter\k)) \to (\bar x, f(\bar x))$  and $\liminf_{\k\to\infty}\, \vnorm[-]{\partial f(x\iter\k)} = 0$, then $0\in \partial f(\bar x)$.
\end{LEM}
Let $f$ be smooth in a neighborhood of $\bar x$ and $g$ be lsc and finite at $\bar x$, then the following holds:
\begin{equation} \label{eq:appdx:sum-formula-subdiff}
  \partial (f+g)(\bar x) = \partial g(\bar x) + \nabla f(\bar x)\,. 
\end{equation}

\section{Accelerated Variants} \label{sec:acc-aFista}

\subsection{Variants with $O(1/\k^2)$-convergence rate} \label{sec:aFista-acc}

In this section, we introduce two variants of our method (with
$\bvec\beta=\beta\in\R$) for convex optimization problems \eqref{eq:gen-opt},
which have a convergence rate of $O(1/\k^2)$. The two methods are variants in
the sense that the update step \eqref{eq:aFista-step-gen-ho} (or an inexact
version) is embedded into another algorithmic strategy. Both cases are related
to the standard FISTA method \cite{BT09b,Tseng08}. 

We define a sequence $\seq[\k\in\N]{\theta\iter\k}$ with $\theta\iter0=\theta\iter{-1} \in (0,1]$ and $\theta\iter\kp\in(0,1]$ such that
\begin{equation} \label{eq:theta-seq-cond}
  \frac{1-\theta\iter\kp}{\theta\iter\kp^2} \leq \frac{1}{\theta\iter\k^2} 
\end{equation}
holds, e.g.  $\theta\iter\k=2/(\k+2)$.

\subsubsection{Monotone FISTA version} \label{subsec:MFZeroSR1}

We use the idea of \cite{LL15} (see also \cite{BT09c}) to obtain an acceleration. 

\mybox{
  \begin{ALG}[Adaptive Monotone FISTA] \ \label{alg:amfista}
    We generate a sequence $\seq[\k\in\N]{z\iter\k}$ starting at some
    $z\iter0\in\R^N$ that obeys the following update scheme: (set $z\iter{-1}=\akp{x}\iter0=z\iter0$)
    \begin{eqnarray}
        \akk{y}\iter\k &=& z\iter\k + \frac{\theta\iter\k(1-\theta\iter\km)}{\theta\iter\km}(z\iter\k-z\iter\km) + \frac{\theta\iter\k}{\theta\iter\km}(\akp{x}\iter\k-z\iter\k) \\
        \akk{x}\iter\kp &=&   \argmin_{x}\ell_f^g(x;\akk{y}\iter\k) + \frac 12 \vnorm{x-\akk{y}\iter\k}^2\\
        y\idx{\beta}\iter\k &=& z\iter\k + \beta(z\iter\k-z\iter\km) \\
        \label{eq:mfista-var-opt-beta-step}
        \akp{x}\iter\kp &=& \argmin_{x} \min_{\beta}\, \ell_f^g(x;y\idx{\beta}\iter\k) + \frac 12 \vnorm{x-y\idx{\beta}\iter\k}^2 \\
        \label{eq:mfista-var-select-step}
        z\iter\kp &=& \begin{cases} \akp{x}\iter\kp\,,&\text{if}\ f^g(\akp{x}\iter\kp) \leq f^g(\akk{x}\iter\kp) \\
                                    \akk{x}\iter\kp\,,&\text{if}\ f^g(\akp{x}\iter\kp) > f^g(\akk{x}\iter\kp)
                      \end{cases}
    \end{eqnarray}
  \end{ALG}
}
\begin{REM}
  \eqref{eq:mfista-var-opt-beta-step} is an update step as in our method \eqref{eq:aFista-step-gen-ho} with $d\iter\k = z\iter\k-z\iter\km$.
\end{REM}
\begin{PROP} \label{prop:mfista-version-rate}
  Let $\seq[\k\in\N]{\theta\iter\k}$ obey \eqref{eq:theta-seq-cond}. The sequence $\seq[\k\in\N]{z\iter\k}$ obeys the following rate of convergence with respect to the objective values: 
  \[
    f^g(z\iter\k) - f^g(\opt{x}) \leq \frac{\theta\iter\k^2}{1-\theta\iter\k} \Bigg(\frac L2 \vnorm{z\iter 0 - \opt{x}}^2 + \frac{1-\theta\iter0}{\theta\iter0} (f^g(z\iter0) - f^g(\opt{x})) \Bigg) \in O(1/\k^2)\,.
  \]
\end{PROP}
The proof is in Section~\ref{proof:prop:mfista-version-rate}.

\subsubsection{Tseng-like version} \label{subsec:TsengZeroSR1}

The following variant does not require a comparison of the objective
value, it relies on a comparison of the value of the proximal linearization.
This algorithm is closely connected to \cite[Algorithm 1]{Tseng08}. Here, we
need an auxiliary sequence $\seq[\k\in\N]{z\iter\k}$ and consider the following
update scheme:
\mybox{
  \begin{ALG}[Adaptive Tseng Acceleration] \ \label{alg:atseng}
    We generate a sequence $\seq[\k\in\N]{\akp{x}\iter\k}$ starting at some
    $\akp{x}\iter0\in\R^N$ that obeys the following update scheme: (set $\akk{x}\iter0=\akp{x}\iter0$)
    \begin{eqnarray}
    \hspace*{-4ex}  \akk{y}\iter\k &=& (1-\theta\iter\k)\akp{x}\iter\k + \theta\iter\k \akk{x}\iter\k \\
      \label{eq:tseng-version-sur-opt}
    \hspace*{-4ex}  \akk{x}\iter\kp &=& \arg\min_{x}\, \lin{f}{g}{x}{\akk{y}\iter\k} + \theta_\k\frac L2 \vnorm{x - \akk{x}\iter\k}^2 \\
      \label{eq:tseng-version-post-comb}
    \hspace*{-4ex}  z\iter\kp &=& (1-\theta\iter\k)\akp{x}\iter\k + \theta\iter\k \akk{x}\iter\kp \\
      \label{eq:tseng-version-selection-step}
    \hspace*{-4ex}  \text{find}\ (\akp{x}\iter\kp,\akp{y}\iter\k) &\st& \lin{f}{g}{\akp{x}\iter\kp}{\akp{y}\iter\k} + \frac L2 \vnorm{\akp{x}\iter\kp - \akp{y}\iter\k}^2 \leq 
      \lin{f}{g}{z\iter\kp}{\akk{y}\iter\k} + \frac L2 \vnorm{z\iter\kp - \akk{y}\iter\k}^2
    \end{eqnarray}
  \end{ALG}
}
\begin{REM}
    Our scheme \eqref{eq:aFista-step-gen-ho} is hidden in \eqref{eq:tseng-version-selection-step}. We can choose any point $\akp{y}\iter\k$. In fact, if we use the exact version of our update scheme \eqref{eq:aFista-step-gen-ho} with $d\iter\k = \akk{x}\iter\k-\akp{x}\iter\k$, there is no need to check this inequality. 
\end{REM}
\begin{PROP} \label{prop:tseng-version-rate}
  Let $\seq[\k\in\N]{\theta\iter\k}$ obey \eqref{eq:theta-seq-cond}. The sequence $\seq[\k\in\N]{\akp{x}\iter\k}$ obeys the following rate of convergence with respect to the objective values:
  \[
    f^g(\akp{x}\iter\k) - f^g(\opt{x}) \leq \frac{\theta\iter\k^2}{1-\theta\iter\k} \Bigg(\frac L2 \vnorm{x\iter0 - \opt{x}}^2 + \frac{1-\theta\iter0}{\theta\iter0} (f^g(\akp{x}\iter0) - f^g(\opt{x})) \Bigg) \in O(1/\k^2)\,.
  \]
\end{PROP}
The proof is in Section~\ref{proof:prop:tseng-version-rate}.

%

\section{Proofs}

\subsection{Proof of Theorem~\ref{thm:aFista-equiv-ho}} \label{proof:thm:aFista-equiv-ho}
  Since $f$ is a quadratic function, we have 
  \[
    f(y\idx{\bvec\beta}) = f(\ak{x}) + \scal{\nabla f(\ak{x})}{y\idx{\bvec\beta}-\ak{x}} + \frac {1}{2}\scal{y\idx{\bvec\beta}-\ak{x}}{\bmat H (y\idx{\bvec\beta}-\ak{x})} \,.
  \]
  and
  \[
    \nabla  f(y\idx{\bvec\beta}) = \nabla f(\ak{x}) + \bmat H(y\idx{\bvec\beta} - \ak{x}) \,.
  \]
  We plug these equations into the objective in \eqref{eq:aFista-step-gen-ho}, use $\bmat D\bvec\beta=y\idx{\bvec\beta}-\ak{x}$, and obtain
  \[
    \begin{split}
    \ell_f^g(x; y\idx{\bvec\beta}) + \frac 1{2} \vnorm[\bmat\Tau]{x-y\idx{\bvec\beta}}^2
    =&\ g(x) + f(\ak{x}) + \scal{\nabla f(\ak{x})}{y\idx{\bvec\beta}-\ak{x}} + \frac 12 \scal{\bmat D \bvec\beta}{\bmat H \bmat D \bvec\beta} \\
    &\ + \scal{\nabla f(\ak{x})}{x-y\idx{\bvec\beta}} + \scal{x-y\idx{\bvec\beta}}{\bmat H\bmat D\bvec\beta}\\
    &\ + \frac 12 \vnorm[\bmat\Tau ]{x-\ak{x}}^2 + \frac 12 \vnorm[\bmat\Tau ]{\bmat D\bvec\beta}^2 - \scal{x-\ak{x}}{\bmat D\bvec\beta}_{\bmat\Tau } \\ 
    =&\ g(x) + f(\ak{x}) + \scal{\nabla f(\ak{x})}{x-\ak{x}} + \frac 12 \vnorm[\bmat\Tau ]{x-\ak{x}}^2 \\
    &\ - \frac 12 \scal{\bmat D \bvec\beta}{\bmat H \bmat D \bvec\beta} + \scal{x-\ak{x}}{\bmat H\bmat D\bvec\beta}   \\
    &\ + \frac 12 \scal{\bmat D\bvec\beta}{\bmat\Tau  \bmat D\bvec\beta} - \scal{x-\ak{x}}{\bmat\Tau \bmat D\bvec\beta} \\
    =&\ g(x) + f(\ak{x}) + \scal{\nabla f(\ak{x})}{x-\ak{x}} + \frac 12 \vnorm[\bmat\Tau ]{x-\ak{x}}^2 \\
    &\ + \frac 12 \scal{\bmat D \bvec\beta}{\bmat M \bmat D \bvec\beta} - \scal{x-\ak{x}}{\bmat M\bmat D\bvec\beta}  \,.
    \end{split}
  \]
  Since $\bmat M\in\spdmat(N)$ and $\bmat D$ has full column rank, we have $\bmat D^\top \bmat M \bmat D\in \spdmat(\nbeta)$. Therefore, this function is strongly convex and differentiable with respect to $\bvec\beta$. The unique optimal choice for $\bvec\beta$, which we denote by $\bvec\beta^*$ can be found by the first order optimality condition:
  \begin{equation} \label{eq:opt-beta-gamma-ho}
    \bmat D^\top \bmat M \bmat D\bvec\beta = \bmat D^\top \bmat M (x-\ak{x}) \,.
  \end{equation}
  Using the optimal $\bvec\beta^*$, the objective in \eqref{eq:aFista-step-gen-ho} reads as follows:
  \[
    \begin{split}
      \ell_f^g(x; y\idx{\opt{\bvec\beta}}) + \frac 1{2} \vnorm[\bmat\Tau]{x-y\idx{\opt{\bvec\beta}}}^2 
      =&\  \ell_f^g(x;\ak{x}) + \frac 12 \vnorm[\bmat\Tau ]{x-\ak{x}}^2 - \frac 12 \scal{\bmat D \bvec\beta^*}{\bmat M \bmat D \bvec\beta^*} \\
      =&\  \ell_f^g(x;\ak{x}) + \frac 12 \vnorm[\bmat Q]{x-\ak{x}}^2 \,,
    \end{split}
  \]
  and the representation in the statement follows directly. The rank of $\bmat U$ is obviously $\nbeta$ and the inversion follows from the rank-$\nbeta$ generalization of the Sherman--Morrison--Woodbury formula.

  It remains to show that indeed $\bmat Q\in \spdmat(N)$. Observe that
  \[
      \bmat Q \bmat D 
    = \bmat T\bmat D -  \bmat U^\top \bmat U \bmat D
    = \bmat T\bmat D -  \bmat M \bmat D 
    = \bmat H \bmat D\quad \text{and}\ \bmat H \in \spdmat(N)
  \]
  holds, which shows that $\scal{x}{\bmat Q x}>0$ for all $x\in \im(\bmat D)$. However, for $x\in \im(\bmat D)^\bot=\ker(\bmat D^\top)$ follows that $\bmat U=0$ and the result follows from $\bmat T\in \spdmat(N)$.
\qed

\subsection{Proof of Remark~\ref{rem:know-H-on-imQ}} \label{proof:rem:know-H-on-imQ}

We define $\DM:=\bmat D^\top \bmat M \bmat D$ and make the following computation for $\bmat Q$:
\[
  \begin{split}
    \bmat Q 
    =&\ \bmat T - \bmat U^\top \bmat U \\
    =&\ \bmat T - \bmat M \bmat D \DM^{-1} (\bmat M \bmat D)^\top \\
    =&\ \bmat T - (\bmat T \bmat D - \bmat Y) \big[\bmat D^\top (\bmat T\bmat D - \bmat Y) \big]^{-1} (\bmat T\bmat D - \bmat Y)^\top  \\
  \end{split}
\]
For $\bmat Q^{-1}$, we observe that $\bmat T^{-1}\bmat M \bmat D=\bmat D - \bmat T^{-1}\bmat Y$ and the following holds:
\[
  \begin{split}
  \matid - \bmat U \bmat\Tau ^{-1} \bmat U^\top 
  =&\  \DM^{-1/2} \big[\DM - \bmat D^\top\bmat M \bmat T^{-1} \bmat M \bmat D \big] \DM^{-1/2} \\
  =&\  \DM^{-1/2} \big[\bmat D^\top(\bmat M - \bmat M \bmat T^{-1} \bmat M) \bmat D \big] \DM^{-1/2} \\
  =&\  \DM^{-1/2} \big[\bmat D^\top(\bmat M - \bmat M \bmat T^{-1} (\bmat T - \bmat H)) \bmat D \big] \DM^{-1/2} \\
  =&\  \DM^{-1/2} \big[\bmat D^\top \bmat M \bmat T^{-1} \bmat H \bmat D \big] \DM^{-1/2} \\
  =&\  \DM^{-1/2} \big[\bmat D^\top (\matid - \bmat H \bmat T^{-1}) \bmat Y \big] \DM^{-1/2} \\
  =&\  \DM^{-1/2} \big[(\bmat D - \bmat T^{-1} \bmat Y)^\top \bmat Y \big] \DM^{-1/2} \\
  \end{split}
\]
Combining these two identities in \eqref{eq:thm:aFista-equiv-ho-Qmatrix-inv} yields the result.
\qed

\subsection{Proof of Proposition~\ref{prop:conv-nc}} \label{proof:prop:conv-nc}

  We start with the following estimation:
  \[
    \begin{split}
      f^g(x\iter\kp)
      \overset{\ii1}{\leq}&\  \ell_f^g(x\iter\kp; y\idx{\bvec\beta\iter\k}) + \frac 12 \norm[\bmat L\piter\k]{x\iter\kp-y\idx{\bvec\beta\iter\k}}^2   \\
      =&\  \ell_f^g(x\iter\kp; y\idx{\bvec\beta\iter\k}) + \frac 12 \norm[\bmat T\piter\k]{x\iter\kp-y\idx{\bvec\beta\iter\k}}^2 + \frac 12 \norm[\bmat L\piter\k - \bmat T\piter\k]{x\iter\kp-y\idx{\bvec\beta\iter\k}}^2 \\
      \overset{\ii2}{\leq}&\  \ell_f^g(x; y\idx{\bvec\beta}) + \frac 12 \norm[\bmat T\piter\k]{x-y\idx{\bvec\beta}}^2 - \frac 12 \norm[\bmat T\piter\k - \bmat L\piter\k]{x\iter\kp-y\idx{\bvec\beta\iter\k}}^2 
    \end{split}
  \]
  where $\ii1$ uses the Descent Lemma (Lemma~\ref{lem:gen-Descent-Lemma}) and $\ii2$ holds for all $x\in \R^N$ and $\bvec\beta\in\R^\nbeta$ due to the update step \eqref{eq:update:alg:afista-nc}. Using the fact that $\frac 12 \norm[\bmat T\piter\k - \bmat L\piter\k]{x\iter\kp-y\idx{\bvec\beta\iter\k}}^2 \geq \frac \lbstep2 \norm{x\iter\kp-y\idx{\bvec\beta\iter\k}}^2$, we obtain 
  \begin{equation} \label{eq:proof:prop:conv-nc:A}
      f^g(x\iter\kp) \leq  \ell_f^g(x; y\idx{\bvec\beta}) + \frac 12 \norm[\bmat T\piter\k]{x-y\idx{\bvec\beta}}^2 - \frac \lbstep2  \norm{x\iter\kp-y\idx{\bvec\beta\iter\k}}^2 \quad \forall x\in \R^N\ \text{and}\ \bvec\beta\in\R^\nbeta\,.
  \end{equation}

  Using $x=x\iter\k$ and $\bvec\beta=0$ in \eqref{eq:proof:prop:conv-nc:A}, we obtain the non-increasingness of $\seq[\k\in\N]{f^g(x\iter\k)}$. By boundedness from below, the sequence $\seq[\k\in\N]{f^g(x\iter\k)}$ must converge. Moreover, we conclude that $x\iter\kp-y\idx{\bvec\beta\iter\k}\to 0$ for $\k\to\infty$.

  Now, let $x^*$ be a limit point of $\seq[\k\in\N]{x\iter\kp}$, i.e., there exist $K\subset\N$ such that $x\iter\kp \to x^*$ for $\k\rto{K}\infty$, hence, $y\idx{\bvec\beta\iter\k}\to x^*$ for $\k\rto{K}\infty$ as well. Consider the limit superior for $\k\rto{K} \infty$ on both sides of the inequality in \eqref{eq:proof:prop:conv-nc:A} with $x=x^*$ and $\bvec\beta=\bvec\beta\piter\k$. Using continuity of $f$ and $\nabla f$, we deduce $\limsup_{\k\rto{K}\infty} g(x\iter\kp) \leq g(x^*)$ and, invoking the lower semi-continuity of $g$, we obtain $g(x\iter\kp) \to g(x^*)$, hence $f^g(x\iter\kp)\to f^g(x^*)$ for $\k\rto{K} \infty$.

  The minimization property of $x\iter\kp$ in \eqref{eq:update:alg:afista-nc} and the property \eqref{eq:appdx:sum-formula-subdiff} yield 
  \[
    0 \in \partial g(x\iter\kp) + \nabla f(y\idx{\bvec\beta\iter\k}) + \bmat T\piter\k (x\iter\kp - y\idx{\bvec\beta\iter\k})\,,
  \]
  which is equivalent to 
  \[
    \nabla f(x\iter\kp) - \nabla f(y\idx{\bvec\beta\iter\k}) - \bmat T\piter\k(x\iter\kp - y\idx{\bvec\beta\iter\k}) \in \partial g(x\iter\kp) + \nabla f(x\iter\kp) = \partial f^g(x\iter\kp) \,.
  \]
  Invoking the Lipschitz continuity of $\nabla f$ (use Assumption~\ref{ass:problem} with Lemma~\ref{lem:gen-Descent-Lemma}), we obtain
  \[
    \norm[-]{\partial f^g(x\iter\kp)} \leq \norm[\bmat L+\bmat T\piter\k]{x\iter\kp - y\idx{\bvec\beta\iter\k}} \,.
  \]
  Finally, boundedness of $\seq[\k\in\N]{\bmat T\piter\k}$ and $x\iter\kp-y\idx{\bvec\beta\iter\k}\to 0$ for $\k\to\infty$ yield $\norm[-]{\partial f^g(x\iter\kp)}\to 0$. The closedness property of the limiting subdifferential together with $f^g(x\iter\kp)\to f^g(x^*)$, for $\k\rto{K} \infty$, allows us to deduce that $x^*$ is a stationary point of $f^g$.
\qed

\subsection{Proof of Proposition~\ref{prop:conv-inc}}  \label{proof:prop:conv-inc}

  Following the same argument as in the proof of Proposition~\ref{prop:conv-nc}, using \eqref{prop:conv-inc:inexact-update} instead of \eqref{eq:update:alg:afista-nc}, we deduce that \eqref{eq:proof:prop:conv-nc:A} holds with $x=x\iter\k$ and $\bvec\beta=0$. Non-increasingness, convergence of $\seq[\k\in\N]{f^g(x\iter\k)}$, and $x\iter\kp-y\idx{\bvec\beta\iter\k}\to 0$ as $\k\to\infty$ follow directly. 
  
  Now, let $x^*$ be a limit point of $\seq[\k\in\N]{x\iter\kp}$ and $K\subset\N$ such that $x\iter\kp\to x^*$ for $\k\rto{K}\infty$. The conditions imply continuity of $f^g$ on $\dom g$. As $\seq[\k\in\N]{f^g(x\iter\k)}$ is non-increasing, $x\iter\k\in\dom g$ for all $\k\in \N$, hence $f^g(x\iter\k) \to f^g(x^*)$ for $\k\rto{K}\infty$. Together with the closedness property of the limiting subdifferential and $\norm[-]{\partial f^g(x\iter\kp)}\to 0$, we deduce that $x^*$ is a critical point of $f^g$.
\qed

\subsection{Proof of Proposition~\ref{prop:mfista-version-rate}} \label{proof:prop:mfista-version-rate}

We make the following estimation:
\[
  \begin{split}
  f^g(z\iter\kp) 
  \overset{\ii1}{\leq}&\ f^g(\akk{x}\iter\kp) \\
  \overset{\ii2}{\leq}&\ \lin{f}{g}{\akk{x}\iter\kp}{\akk{y}\iter\k} + \frac L2 \vnorm{\akk{x}\iter\kp - \akk{y}\iter\k}^2 \\
  \overset{\ii3}{\leq}&\ \lin{f}{g}{y}{\akk{y}\iter\k} + \frac L2 \vnorm{y - \akk{y}\iter\k}^2 - \frac L2 \vnorm{y - \akk{x}\iter\kp}^2 \\
  \overset{\ii4}{\leq}&\ \lin{f}{g}{(1-\theta_\k)z\iter\k + \theta_\k x}{\akk{y}\iter\k} + \frac L2 \vnorm{(1-\theta_\k)z\iter\k + \theta_\k x - \akk{y}\iter\k}^2 - \frac L2 \vnorm{(1-\theta_\k)z\iter\k + \theta_\k x - \akk{x}\iter\kp}^2 \\
  \overset{\ii5}{\leq}&\ (1-\theta_\k)\lin{f}{g}{z\iter\k}{\akk{y}\iter\k} + \theta_\k \lin{f}{g}{x}{\akk{y}\iter\k} + \theta_\k^2 \frac L2 \vnorm{(1-\theta_\k)/\theta_\k z\iter\k + x - \akk{y}\iter\k/\theta_\k}^2\\
  &\ - \theta_\k^2\frac L2 \vnorm{(1-\theta_\k)/\theta_\k z\iter\k + x - \akk{x}\iter\kp/\theta_\k}^2 \\
  \overset{\ii6}{\leq}&\ (1-\theta_\k)f^g(z\iter\k) + \theta_\k f^g(x) + \theta_\k^2 \frac L2 (\vnorm{U\iter\k(x)}^2 - \vnorm{U\iter\kp(x)}^2) 
  \end{split}
\]
where $\ii1$ uses \eqref{eq:mfista-var-select-step}, $\ii2$ uses the quadratic (Lipschitz) upper bound, $\ii3$ holds for all $y$ by using $L$-strong convexity of $x\mapsto \lin{f}{g}{x}{\akk{y}\iter\k} + \frac L2 \vnorm{x - \akk{y}\iter\k}^2$ and optimality of $\akk{x}\iter\kp$, $\ii4$ holds for all $x$ by the change of variables $y=(1-\theta_\k)z\iter\k + \theta_\k x$, $\ii5$ holds by convexity of $\linfun fg$ in the first argument and a simple algebraic manipulation, and $\ii6$ holds by defining $U\iter\kp(x):=(1-\theta_\k)/\theta_\k z\iter\k + x - \akk{x}\iter\kp/\theta_\k$ and the definition of\footnote{This equality was actually used to define $\akk{y}\iter\k$.} $\akk{y}\iter\k$.

Rearranging this inequality and setting $x=\opt{x}$ yields
\[
  \frac{1}{\theta\iter\k^2}(f^g(z\iter\kp) - f^g(\opt{x}) - \frac{(1-\theta\iter\k)}{\theta\iter\k^2}(f^g(z\iter\k)-f^g(\opt{x})) \leq \frac L2 (\vnorm{U\iter\k(\opt{x})}^2 - \vnorm{U\iter\kp(\opt{x})}^2) \,,
\]
which by standard arguments shows the result.
\qed

\subsection{Proof of Proposition~\ref{prop:tseng-version-rate}} \label{proof:prop:tseng-version-rate}

  We make the following estimation:
\[
  \begin{split}
  f^g(\akp{x}\iter\kp) 
  \overset{\ii1}{\leq}&\ \lin{f}{g}{\akp{x}\iter\kp}{\akp{y}\iter\k} + \frac L2 \vnorm{\akp{x}\iter\kp - \akp{y}\iter\k}^2 \\
  \overset{\ii2}{\leq}&\ \lin{f}{g}{z\iter\kp}{\akk{y}\iter\k} + \frac L2 \vnorm{z\iter\kp - \akk{y}\iter\k}^2 \\
  \overset{\ii3}{\leq}&\ \lin{f}{g}{(1-\theta\iter\k)\akp{x}\iter\k + \theta\iter\k \akk{x}\iter\kp}{\akk{y}\iter\k} + \frac L2 \vnorm{(1-\theta\iter\k)\akp{x}\iter\k + \theta\iter\k \akk{x}\iter\kp - \akk{y}\iter\k}^2 \\
  \overset{\ii4}{\leq}&\ (1-\theta\iter\k)\lin{f}{g}{\akp{x}\iter\k}{\akk{y}\iter\k}  + \theta\iter\k (\lin{f}{g}{\akk{x}\iter\kp}{\akk{y}\iter\k} + \theta\iter\k\frac L2 \vnorm{\akk{x}\iter\kp - \akk{x}\iter\k}^2) \\
  \overset{\ii5}{\leq}&\ (1-\theta\iter\k)\lin{f}{g}{\akp{x}\iter\k}{\akk{y}\iter\k} + \theta\iter\k (\lin{f}{g}{x}{\akk{y}\iter\k} + \theta\iter\k\frac L2 \vnorm{x - \akk{x}\iter\k}^2 - \theta\iter\k\frac L2 \vnorm{x-\akk{x}\iter\kp}^2) \\
  \overset{\ii6}{\leq}&\ (1-\theta\iter\k)f^g(\akp{x}\iter\k) + \theta\iter\k f^g(x) + \theta\iter\k^2\frac L2( \vnorm{U\iter\k(x)}^2 - \vnorm{U\iter\kp(x)}^2) \\
  \end{split}
\]
where $\ii1$ uses the quadratic (Lipschitz) upper bound, $\ii2$ uses \eqref{eq:tseng-version-selection-step}, $\ii3$ uses \eqref{eq:tseng-version-post-comb}, $\ii4$ uses convexity of $\ell_f^g$ and the definition of\footnote{Again, $\akk{y}\iter\k$ is actually chosen such that this equality holds.} $\akk{y}\iter\k$ , $\ii5$ holds for any $x$ thanks to the optimality of $\akk{x}\iter\kp$ in \eqref{eq:tseng-version-sur-opt} and the $L$-strong convexity of $x\mapsto \lin{f}{g}{x}{\akk{y}\iter\k} + \theta\iter\k\frac L2 \vnorm{x - \akk{x}\iter\k}^2$, $\ii6$ holds by defining $U\iter\k(x):= x - \akk{x}\iter\k$  and by convexity of $\ell_f^g$. The statement follows analogously to Proposition~\ref{prop:mfista-version-rate}.
\qed

{\small
\bibliographystyle{ieee}
\bibliography{ochs}
}

\end{document}